\documentclass[12pt]{article}
\usepackage[latin1]{inputenc}
\usepackage{amssymb,amsmath,amsfonts}
\usepackage[notcite,notref]{showkeys}
\newtheorem{theorem}{Theorem}

\newtheorem{lemma}[theorem]{Lemma}

\newtheorem{remarks}[theorem]{Remarks}

\newcommand{\R}{\mathbb{R}}
\newcommand{\Q}{\mathbb{Q}}
\newcommand{\Sf}{\mathbb{S}}
\newcommand{\C}{\mathbb{C}}

\newcommand{\spa}{\mbox{span}}

\newcommand{\po}{{\hspace*{-1ex}}{\bf .  }}

\newcommand{\E}{{\cal E}}
\def\d{{\partial}}
\def\<{{\langle}}
\def\>{{\rangle}}

\def\d{\partial}

\def\a{\alpha}

\def\bea{\begin{eqnarray*} }
\def\eea{\end{eqnarray*} }
\def\be{\begin{equation} }
\def\ee{\end{equation} }

\def\proof{\noindent\emph{ Proof: }}
\def\qed{\ifhmode\unskip\nobreak\fi\ifmmode\ifinner
\else\hskip5 pt \fi\fi\hbox{\hskip5 pt \vrule width4 pt
height6 pt  depth1.5 pt \hskip 1pt }}

\begin{document}

\title{A representation for pseudoholomorphic\\ surfaces in spheres}
\author{M. Dajczer and Th. Vlachos}
\date{}
\maketitle

\begin{abstract} We give a local representation for the  pseudoholomorphic
surfaces in Euclidean spheres in terms of holomorphic data.  Similar to the 
case of the generalized Weierstrass representation of Hoffman and 
Osserman, we assign such a surface in $\Sf^{2n}$ to a given set of $n$ 
holomorphic functions defined on a simply-connected domain in $\C$. 
\end{abstract}

A basic tool in the study of minimal surfaces in Euclidean space  with codimension 
higher than one is the generalized Weierstrass representation 
introduced by  Hoffman and Osserman in \cite{HF}. Roughly speaking, to a given 
set of $n$ holomorphic functions defined on a simply-connected domain in $\C$
it assigns a generalized minimal surface in $\R^{n+1}$. 
Conversely,  any  simply-connected minimal surface in $\R^{n+1}$ 
can be obtained in this way.

For minimal surfaces in Euclidean spheres, there is no  representation
similar to the one in Euclidean space and, maybe, this is not even possible. 
But inspired by results and methods  in \cite{CCh}, \cite{df1} and \cite{DG}, 
we provide such a representation for a class of minimal surfaces, namely, 
the generalized pseudoholomorphic 
surfaces in Euclidean spheres. As above, to a set of $n$ holomorphic 
functions in a simply-connected domain in $\C$ and some constants 
of integration we assign a generalized pseudoholomorphic  surfaces in 
$\Sf^{2n}$. Here and above, by the term generalized we mean that the 
metric induced on the surface is singular at most at isolated points.

In a seminal paper due to Calabi \cite{Cal} (see also Barbosa \cite{Ba}) it 
was shown that any two-dimensional sphere minimally immersed in an Euclidean 
sphere lays in even substantial codimension and has to be pseudoholomorphic. 
According to Calabi's definition, that a minimal surface is pseudoholomorphic 
means that when parametrized by an isothermal coordinate $z$
the surface $f\colon\Sf^2\to\Sf^{2n}\subset\R^{2n+1}$  satisfies 
$$
\<\d^jf,\d^kf\>=0,\;\;j+k>0,
$$
with respect to the symmetric product where $\d=\partial/\partial z$ and 
$\d^0f=f$. 

There are other definitions for minimal surfaces in  spheres to be pseudoholomorphic
that are equivalent to the above.
For instance, one may require the ellipses of curvature of any order and at any point 
to be circles; see part $(iii)$ in Remarks \ref{remarks}. 
In turn, this is equivalent to ask the Hopf differentials
of all orders to  vanish identically; see \cite{Vl} for details.

We point out that even in the compact case  there are plenty of pseudoholomorphic 
surfaces other than topological spheres.  For instance, this is the case of many 
pseudoholomorphic surfaces in the nearly Kaehler sphere $\Sf^6$.  
These surfaces were introduced by Bryant \cite{Br} and have been intensively studied.
\vspace{1,5ex}

In the sequel, we first describe an algorithm that leads to the representation of 
the pseudoholomorphic surfaces and after that we state our main result.
\vspace{1,5ex}

 Let $\beta_s\colon U\to\C$, $0\leq s\leq n-1$, be nonzero holomorphic functions 
defined on a simply-connected domain $U\subset \C$. Setting $\a_0=\beta_0$,  define 
subsequently isotropic maps with respect to the symmetric product between complex 
vectors  
$\a_r\colon\, U\to\C^{2r+1}$, $0\leq r\leq n-1$,  by
$$
\a_{r+1}=\beta_{r+1}\left(1-\phi_r^2,i(1+\phi_r^2),2\phi_r\right)
$$
where $\phi_r=\int^z\a_rdz$ and $\beta_n=1$. 
Now define $F_s\colon U\to\C^{2n+1},\;\;1\leq s\leq n+1$, by
\be\label{F1}
F_1=\a_n\;\;\;\mbox{and}\;\;\;
F_{r+1}=\d F_r-\dfrac{1}{\|F_r\|^2}\<\d F_r,\bar{F}_r\>F_r,\;1\leq r\leq n.
\ee
Finally, let $g\colon L^2\to\Sf^{2n}$ be the map given by
$$
g=\frac{1}{\|{\mathrm{Re}}(F_{n+1})\|}{\mathrm {Re}}(F_{n+1})
$$
where $L^2=U$ with the induced metric.

\begin{theorem}\po\label{main}
The map $g\colon L^2\to\Sf^{2n}$ parametrizes a generalized pseudoholomorphic 
\mbox{surface}. Conversely, any pseudoholomorphic surface in $\Sf^{2n}$ can 
be locally obtained in this way.
\end{theorem}

To conclude the paper, we show how the above parametrization can be used
to  parametrize the real Kaehler hypersurfaces in Euclidean space as well
as a class of ruled minimal submanifolds with codimension two in spheres.
\newpage

\section{The proof}

We first collect some basic facts and definitions about minimal surfaces in space
forms and refer to \cite{df1} for more details.
\vspace{1,5ex}

Let $f\colon L^2\to\Q^N$ denote an isometric immersion of a
two-dimensional  Riemannian manifold into an $N$-dimensional  space form.
The $k^{th}$\emph{-normal space} of $f$ at $x\in L^2$ for $k\geq 1$ is
defined as
$$
N^f_k(x)=\spa\{\a_f^{k+1}(X_1,\ldots,X_{k+1}):X_1,\ldots,X_{k+1}\in T_xL\}
$$
where
$$
\a_f^s\colon TL\times\cdots\times TL\to N_fL,\;\; s\geq 3,
$$
denotes the symmetric tensor called the $s^{th}$\emph{-fundamental form} given
inductively by
$$
\a_f^s(X_1,\ldots,X_s)=\left(\nabla^\perp_{X_s}\ldots
\nabla^\perp_{X_3}\a_f(X_2,X_1)\right)^\perp
$$
and $\a_f\colon TL\times TL\to N_fL$ stands for the standard second fundamental
form of $f$ with values in the normal bundle.
Here  $\nabla^{\perp}$ denotes the induced connection in the normal bundle $N_fL$ of $f$
and $(\;\;)^\perp$ means taking the projection onto the normal complement of
$N^f_1\oplus\ldots\oplus N^f_{s-2}$ in $N_fL$.

A surface $f\colon L^2\to\Q^N$ is called \emph{regular} if for each $k$ the 
subspaces $N^f_k$ have constant dimension and thus form normal subbundles
of the normal bundle.
Notice that regularity is always verified  along connected components of 
an open dense subset of $L^2$.

Assume that $f\colon L^2\to\Q^N$ is a substantial regular minimal surface.
Substantial means that the codimension cannot be reduced.
Then, the normal bundle of $f$ splits as
$$
N_fL=N_1^f\oplus N_2^f\oplus\dots\oplus N_m^f,\;\;\; m=[(N-1)/2],
$$
since all higher normal bundles  have rank two except possible the last
one that has rank one if $N$ is odd.
The \emph{$k^{th}$-order curvature ellipse}
$\E^f_k(x)\subset N^f_{k}(x)$ at $x\in L^2$  for each
$N_k^f$, $1\leq k\leq m$, is defined  by
$$
\E^f_k(x) = \{\a_f^{k+1}(Z^{\varphi},\ldots,Z^{\varphi})\colon\,
Z^{\varphi}=\cos\varphi Z+\sin\varphi JZ\;\mbox{and}\;\varphi\in\Sf^1\}
$$
where $Z\in T_xL$ is any vector of unit length and $J$ is the complex structure
in $T_xL$.

A substantial regular minimal surface $f\colon L^2\to\R^{2n+1}$ is called 
\emph{isotropic} if at any point the ellipses of curvature $\E^f_k(x)$,
$1\leq k\leq n-1$, are circles. From the results in \cite{CCh} and \cite{DG}
we have that such a surface can be locally parametrized as 
$f=\mathrm{Re}(\psi\phi_n)$ where $\psi$ is any nowhere vanishing holomorphic 
function.  The condition that the ellipses
of curvature are circles is equivalent to the vector fields 
$F_1,\d F_1,\dots,\d^{n-1}F_1$ being orthogonal and isotropic.
Notice that different functions $\psi$ yield isotropic surfaces
with the same Gauss map.
\vspace{1,5ex}
\newpage 

For the proof of Theorem \ref{main} we first give several lemmas.

\begin{lemma}\po\label{bas}  The following facts hold:
\begin{itemize}
\item[(i)]  Outside the zeros $F_1,\dots,F_n$  span a maximal 
isotropic subspace of $\C^{2n+1}$.
\item[(ii)] The vectors $F_1,\dots,F_{n+1}$ are orthogonal with 
respect to the Hermitian product.
\item[(iii)] The vectors $F_{n+1}, \bar{F}_{n+1}$ are collinear. 
\item[(iv)] It holds that
\be\label{F}
F_{s+1}=\d F_s-\d(\log\|F_s\|^2)F_s,\;\; 1\leq s\leq n.
\ee
\item[(v)] It holds that
\be\label{Fbar}
\d\bar{F}_s=-\frac{\|F_s\|^2}{\|F_{s-1}\|^2}\bar{F}_{s-1},\;\; 2\leq s\leq n.
\ee
\item[(vi)]  The zeros of $F_1,\dots,F_{n+1}$ are isolated.
\end{itemize}
\end{lemma}

\proof 
\noindent $(i)$ This follows from
$$
\spa_{\,\C} \{F_1,\dots, F_n\}=\spa_{\,\C} \{F_1,\d F_1,\dots,
\d^{n-1}F_1\}.
$$
and that $n$ is the dimension of any maximal isotropic subspaces in $\C^{2n+1}$.

\noindent $(ii)$  We claim that
\be\label{fund}
\d\bar{F}_s \in \spa_{\,\C} \{\bar{F}_1,\dots, \bar{F}_{s-1}\},\;\; 2\leq s\leq n.
\ee
In fact, since $F_1$ is holomorphic we see that
$$
\d\bar{F}_2=-\d\left(\<\bar\d\bar{F}_1,F_1\>/\|F_1\|^2\right)\bar{F}_1.
$$
Assume that the claim holds for any $2\leq s\leq k$.  Then, we have
$$
\d\bar{F}_{k+1}=\d\bar\d \bar{F}_k-\d\left(\<\bar\d\bar{F}_k,F_k\>/\|F_k\|^2\right)
\bar{F}_k-\left(\<\bar\d\bar{F}_k,F_k\>/\|F_k\|^2\right)\d\bar{F}_k.
$$
By assumption, we can write 
$$
\d\bar{F}_k=\lambda_1\bar{F}_1+\cdots+\lambda_{k-1}\bar{F}_{k-1}.
$$
Thus,
$$
\d\bar\d\bar{F}_k=\sum_{j=1}^{k-1}\bar\d(\lambda_j)\bar{F}_j
+\sum_{j=1}^{k-1}\lambda_j\bar\d\bar{F}_j
$$
and since  $\bar\d\bar{F}_j\in\spa_{\,\C}\{\bar{F}_j,\bar{F}_{j+1}\}$, 
the claim follows.

That $\<F_2,\bar{F}_1\>=0$ is immediate.
Assume that $F_1,\dots,F_s$ are orthogonal  with respect to the Hermitian product. 
For $1\leq j\leq s-1$  we obtain using (\ref{fund}) that
\bea
\<F_{s+1},\bar{F}_j\>\!\!\!&=&\!\!\!\<\d F_s,\bar{F}_j\>
-\dfrac{1}{\|F_s\|^2}\<\d F_s, \bar{F}_j\>\<F_s,\bar{F}_j\>\\
\!\!\!&=&\!\!\! \d\< F_s, \bar{F}_j\>-\< F_s, \d\bar{F}_j\>\\
\!\!\!&=&\!\!\!0.
\eea
Since also $\<F_{s+1},\bar{F}_s\>=0$ is immediate, the result follows. 
\vspace{1ex}

\noindent $(iii)$ It suffices to show that $F_{n+1}$ is orthogonal to $F_j,\bar{F}_j$, 
$1\leq j\leq n$,  with respect to the Hermitian inner product.  
Then the same holds for $\bar{F}_{n+1}$ and the result follows.
We have from part $(ii)$ that $F_{n+1}$ is orthogonal to $F_j$, $1\leq j\leq n$.
The orthogonality with respect to $\bar{F}_j$, $1\leq j\leq n$, follows from $(i)$.
\vspace{1ex}

\noindent $(vi)$ Immediate using (\ref{fund}).
\vspace{1ex}

\noindent $(v)$   Follows easily using part $(ii)$ and (\ref{fund}).
\vspace{1ex}

\noindent $(v)$ Since $F_1$ is holomorphic, we have using (\ref{Fbar}) that
$$
\bar\d\left(F_1\wedge\cdots\wedge F_{n+1}\right)=0,
$$
and the result follows.\qed

\begin{lemma}\po\label{bundles}
The complexified tangent and normal bundles of 
$f=\mathrm{Re}(\phi_n)\colon U\to\R^{2n+1}$  are given by
\bea
f_*TU\otimes\C\!\!\!&=&\!\!\!\spa_{\,\C}\{F_1,\bar{F}_1\},\\
N_s^f\otimes\C\!\!\!&=&\!\!\!\spa_{\,\C}\{F_{s+1},\bar{F}_{s+1}\},
\;\; 1\leq s\leq n-1,\\
N_n^f\otimes\C\!\!\!&=&\!\!\!\spa_{\,\C} \{F_{n+1}\}.
\eea
\end{lemma}

\proof The first equality follows from 
$$
2\d f=\d\phi_n=F_1.
$$ 
We claim that the higher fundamental forms satisfy
\be\label{a}
2\a_f^{s+1}(\d,\dots,\d)=F_{s+1},\;\; 1\leq s \leq n.
\ee
We proceed by induction. We have that
$$
2\a_f(\d,\d)=2(\d\d f)^\perp=(\d F_1)^\perp=\d F_1-(\d F_1)^{f_*TU}.
$$
Computing $(\d F_1)^{f_*TU}$ in the real tangent base 
$\{F_1+\bar{F}_1, i(F_1-\bar{F}_1)\}$ 
gives
$$
(\d F_1)^{f_*TU}=\dfrac{1}{\|F_1\|^2}\<\d F_1, \bar{F}_1\>F_1,
$$
and (\ref{a}) follows for $s=1$.

Now assume by induction that (\ref{a}) holds for any $1\leq k\leq s$. 
A similar computation as above  using parts $(i)$ and $(ii)$ 
of Lemma \ref{bas} gives
\bea
2\a_f^{s+2}(\d,\dots,\d)\!\!\!&=&\!\!\!\d F_{s+1}
-(\d F_{s+1})^{f_*TU}-\sum_{j=1}^s(\d F_{s+1})^{N^f_j}\\
\!\!\!&=&\!\!\!\d F_{s+1}- (\d F_{s+1})^{N^f_s}\\
\!\!\!&=&\!\!\!\d F_{s+1}- \dfrac{\<\d F_{s+1}, 
\bar{F}_{s+1}\>}{\|F_{s+1}\|^2}F_{s+1}\\
\!\!\!&=&\!\!\!F_{s+2},
\eea
and this completes the proof.\qed

\begin{lemma}\po\label{polar}
The  surface $g\colon L^2\to\Sf^{2n}$ is minimal and we have:
\begin{eqnarray}
g_*TL\otimes\C\!\!\!&=&\!\!\!\spa_{\,\C} \{F_n,\bar{F}_n\},\label{first}\\
N_s^g\otimes\C\!\!\!&=&\!\!\!\spa_{\,\C} \{F_{n-s},\bar
F_{n-s}\},\;\; 1\leq s \leq n-1,\label{second}
\end{eqnarray}
and
\be\label{ap}
\a_g^{s+1}(\d,\dots,\d)
=\frac{(-1)^{s+1}}{\|F_{n-s}\|^2}\<g,F_{n+1}\>\bar{F}_{n-s},\;\; 1\leq s \leq n-1.
\ee
\end{lemma}

\proof We easily obtain using Lemma \ref{bas} that $\d F_{n+1},\d\bar{F}_{n+1}$ 
are orthogonal to $F_s, \bar{F}_s$, $1\leq s \leq n-1$.  From part  $(iii)$ of 
Lemma \ref{bas} we obtain that $\d F_{n+1},\d\bar{F}_{n+1}$ 
are orthogonal to $F_{n+1}$, and this gives (\ref{first}). 

Using (\ref{first}) and Lemma \ref{bas} it follows easily that
\be\label{tp}
\d g=-\frac{\<g,F_{n+1}\>}{\|F_n\|^2}\bar{F}_n.
\ee
Hence,
$$
\a_g(\d,\bar\d)=(\d\bar\d g)^\perp
=-\frac{\<g,\bar{F}_{n+1} \>}{\|F_n\|^2}(\d F_n)^\perp.
$$
However, 
\bea
(\d F_n)^\perp\!\!\!&=&\!\!\!\d F_n-\<\d F_n,g\>g-(\d F_n)^{g_*TL}\\
\!\!\!&=&\!\!\!\d F_n-\<F_{n+1},g\>g-\d(\log \|F_n\|^2)F_n\\
\!\!\!&=&\!\!\!F_{n+1}-\<F_{n+1},g\>g\\
\!\!\!&=&\!\!\!0,
\eea
and thus $g$ is minimal.

To prove (\ref{second}) and (\ref{ap}) we proceed by induction. 
Using  Lemma \ref{bas} and (\ref{tp}), we obtain
\bea
\a_g(\d,\d)
\!\!\!&=&\!\!\!\d\d g-\<\d\d g,g\>g-(\d\d g)^{g_*TL}\\
\!\!\!&=&\!\!\!-\d\left(\frac{\<g,F_{n+1}\>}{\|F_n\|^2}\right)\bar{F}_n+\frac{\<g,F_{n+1}
\>}{\|F_{n-1}\|^2}\bar{F}_{n-1}\\
\!\!\!&&\!\!\!+\d\left(\frac{\<g,F_{n+1} \>}{2\|F_n\|^2}\right)(F_n+\bar{F}_n)
-\d\left(\frac{\<g,F_{n+1}\>}{2\|F_n\|^2}\right)(F_n-\bar{F}_n)\\
\!\!\!&=&\!\!\!\frac{1}{\|F_{n-1}\|^2}\<g,F_{n+1}\>\bar{F}_{n-1}.
\eea
Now assume that (\ref{ap}) holds for any $1\leq k\leq s$. 
Similarly, we have
\bea
\a_g^{s+2}(\d,\dots,\d)\!\!\!&=&\!\!\!\d\a_g^{s+1}(\d,\dots,\d) 
-(\d\a_g^{s+1}(\d,\dots,\d))^{g_*TL}-\sum_{j=1}^s(\d\a_g^{s+1}(\d,\dots,\d))^{N^g_j}\\
\!\!\!&=&\!\!\!(-1)^{s+1}\d\left(\frac{\<g,F_{n+1} \>}{\|F_{n-s}\|^2}\right)\bar{F}_{n-s}
+\frac{(-1)^s}{\|F_{n-s-1}\|^2}\<g,F_{n+1}\>\bar{F}_{n-s-1}\\
\!\!\!&&\!\!\!+(-1)^s\d\left(\frac{\<g,F_{n+1}\>}{2\|F_{n-s}\|^2}\right)(F_{n-s}
+\bar{F}_{n-s})\\
\!\!\!&&\!\!\!+(-1)^{s+1}\d\left(\frac{\<g,F_{n+1}\>}
{2\|F_{n-s}\|^2}\right)(F_{n-s}-\bar{F}_{n-s})\\
\!\!\!&=&\!\!\!\frac{(-1)^s}{\|F_{n-s-1}\|^2}\<g,F_{n+1}\>\bar{F}_{n-s-1},
\eea
and this completes the proof.
\qed\vspace{1.5ex}

\noindent\emph{Proof of Theorem \ref{main}:} The direct statement
follows from the previous lemmas.
For the converse, let $g\colon L^2\to\Sf^{2n}$ be a simply connected 
pseudoholomorphic surface and $z$ a local complex chart. Define $G_0=g$ and
$$
G_{s+1}=\d G_s-\dfrac{\<\d G_s, \bar{G}_s\>}{\|G_s\|^2}G_s,\;\; s \geq 0.
$$
We have that
$$
g_*TL\otimes\C={\spa}_{\C}\{G_1,\bar{G}_1\}.
$$
Clearly, the higher fundamental forms are given by
\be\label{name}
\a_g^{s+1}(\d,\dots,\d)=G_{s+1},\;\; 1\leq s\leq n.
\ee
Thus $G_1,\dots,G_n$ span an isotropic subspace of $\C^{2n+1}$. 
Moreover,
\be\label{normals}
N_s^g\otimes\C=\spa _{\C} \{G_{s+1},\bar{G}_{s+1}\},\;\; 1\leq s \leq n-1.
\ee

We claim that
\be\label{fundG}
\d\bar{G}_s \in {\spa}_{\C} \{\bar{G}_0,\dots, \bar{G}_{s-1}\},\;\; s\geq 1.
\ee
Since $\|g\|=1$ we have  $G_1=\d g$, and hence $\bar G_1=\bar\d g$. Therefore,
$$\vspace{1,5ex}
\d\bar G_1=\d\bar\d g \in\spa_\C\{g\}.
$$
Assume that (\ref{fundG}) holds for $1\leq k\leq s-1$. Then,
$$
\d\bar G_{k+1}=\bar\d\d\bar G_k
-\d\left(\frac{\<\bar\d\bar G_k,G_k\>}{\|G_k\|^2}\right)\bar G_k
-\frac{\<\bar\d\bar G_k,G_k\>}{\|G_k\|^2}\d\bar G_k, 
$$
and the claim follows.

Using (\ref{normals}) we see that $G_1,\dots,G_n$ are orthogonal with respect 
to the  Hermitian product and thus span a maximal isotropic subspace of $\C^{2n+1}$. 
Hence (\ref{fundG}) gives
\be\label{Gbar}\vspace{1,5ex}
\d\bar{G}_s=-\frac{\|G_s\|^2}{\|G_{s-1}\|^2}\bar{G}_{s-1},\;\; s\geq 1.
\ee
We have that
\be\label{vectors}
\C^{2n+1}=\spa_\C\{g,G_1,\ldots,G_n,\bar G_1\ldots,\bar G_n\}
\ee
where all vectors are orthogonal with respect to the  Hermitian product.
It is easy to check that $G_{n+1}$ is orthogonal to all vectors in
(\ref{vectors}) and thus $G_{n+1}=0$ or, equivalently, it holds that
$$
\d G_n=\d \log\|G_n\|^2G_n,
$$
and therefore $\xi=\bar{G}_n/\|G_n\|^2$ is holomorphic.

Consider the map $f\colon L^2\to\R^{2n+1}$ given by
\be\label{f}
f=\mbox{Re}\int\xi dz.\vspace{1,5ex}
\ee
Then,
$$
f_*TL\otimes\C={\spa}_{\C}\{G_n,\bar{G}_n\}.
$$
We have that $f$ is minimal since
$$
\a_f(\d,\bar\d)=(\bar\d\d f)^\perp=0.
$$
From (\ref{Gbar}) we see that
$$
\d^s\bar{G}_n\in{\spa}_{\C}\{\bar{G}_{n-s},\dots,\bar{G}_{n-1}\},
\;\;1\leq s\leq n-1.
$$
Therefore, the subspace 
$$
\spa_\C\{\d\bar{G}_n,\dots,\d^{n-1}\bar{G}_n\}
$$
is isotropic, and hence $f$ is an isotropic surface. 

We have that 
$$
F_1=\d f=\bar{G}_n/\|G_n\|^2.
$$
Using (\ref{Gbar}) we obtain that $F_s$, $1\leq s\leq n$, defined 
by (\ref{F1}) satisfies
$$
F_s\in\spa_\C\{\bar G_n,\ldots,\bar G_{n-s+1}\}.
$$
Thus,
$$
\spa_{\C}\{F_1,\dots,F_n,\bar{F}_1,\dots,\bar{F}_n\}
=\spa_{\C} \{ G_1,\dots,G_n,\bar G_1,\dots,\bar{G}_n\}.
$$
Hence,
\be\label{ex}
g=\frac{1}{\|\mathrm{Re}(F_{n+1})\|}\mathrm{Re}(F_{n+1}),
\ee
and this completes the proof.\qed

\begin{remarks}\po\label{remarks} \emph{$(i)$  Observe that in (\ref{f}) we can 
replace $\xi$ by $\psi\xi$ where $\psi$ is any nowhere
vanishing holomorphic function.  Hence, different holomorphic data 
may generate the same pseudoholomorphic spherical 
surface.\vspace{1ex}\\ 
\noindent $(ii)$ That the first two definitions of pseudoholomorphicity 
given in the introduction are equivalent follows using (\ref{name}).
\vspace{1ex}\\ 
\noindent $(iii)$ The proof of the converse of the theorem can also
be obtained from \cite{df1} by means of the quite different techniques
developed there for the more general context of elliptic surfaces which,  
in particular, do not yield (\ref{ex}).
}\end{remarks}

\section{Applications}

In this section, we first provide a local Weierstrass type representation 
for the real Kaehler  Euclidean hypersurfaces free of flat points. 
\vspace{1,5ex}

The following result was obtained  in \cite{dg1}.
\vspace{1,5ex}

\noindent{\bf Theorem.} \emph{
Let $g\colon L^2 \to \Sf^{2n}$, $n\geq 2$, be a pseudoholomorphic surface 
and $\gamma\in C^{\infty}(L)$  an arbitrary function. 
Then the  induced metric on the open subset of regular 
points of the map $\Psi\colon N_gL\to\R^{2n+1}$ given by
$$
\Psi(x,w) = \gamma (x)g(x) + g_*\nabla\gamma(x)+ w
$$
is Kaehler.
Conversely, any real Kaehler hypersurface 
$f\colon M^{2n} \to \R^{2n+1}, n\geq 2$, free of flat points 
can be locally parametrized in this way.}
\vspace{1,5ex}

We have that
$$
\nabla\gamma=\frac{1}{\|\d\|^2}(\gamma_{\bar z}\d+\gamma_z \bar\d),
$$
$$
g=\frac{1}{\|{\mathrm {Re}}(F_{n+1})\|}{\mathrm {Re}}(F_{n+1})
\;\;\;\mbox{and}\;\;\;
\d g=-\frac{\<g,F_{n+1}\>}{\|F_n\|^2}\bar{F}_n.
$$
Using (\ref{second}) we see that any $w\in N_gL$ can be written as
$
w=\sum_{j=1}^{n-1}{\mathrm {Re}}(w_jF_j)
$
where the $w_j=u_j+iv_j$, $1\leq j\leq n-1$, are complex parameters.

\begin{theorem}\po Any real Kaehler hypersurface 
$f\colon M^{2n} \to \R^{2n+1}, n\geq 2$, without flat points can be 
locally parametrized as
\bea
\Psi(z,w_j) \!\!\!&=&\!\!\! 
\frac{\gamma}{\|{\mathrm{Re}}(F_{n+1})\|}{\mathrm{Re}}(F_{n+1})
-\frac{2}{\|\d\|^2\|F_n\|^2\|{\mathrm{Re}}(F_{n+1})\|}
{\mathrm{Re}}\left({\gamma_z\<{\mathrm {Re}}(F_{n+1}),\bar{F}_{n+1}\>}
F_n\right)\\
\!\!\!&&\!\!\!+\sum_{j=1}^{n-1}(u_j\mathrm{Re}(F_j)-v_j\mathrm {Im}(F_j)).
\eea
\end{theorem}

\proof A straightforward computation using (\ref{tp}).\qed
\vspace{1,5ex}

Next, we describe how to parametrize the minimal ruled submanifolds 
in spheres with codimension two  given in \cite{dv2} but only
when associated to a pseudoholomorphic surface.
\vspace{1,5ex}

Let  $g\colon L^2\to\Sf^{2n}$, $n\geq 3$,  be a pseudoholomorphic surface.
We constructed in \cite{dv2} an associated ruled minimal submanifold 
$F_g\colon M^{2n-2}\to\Sf^{2n}$ by attaching at each point of $g$ 
the totally geodesic 
$(2n-4)$-sphere of $\Sf^{2n}$ whose tangent space at that point is the 
fiber  of the vector bundle  
$\Lambda_g=(N_1^g)^\perp$, that is,
$$
(p,w)\in\Lambda_g\mapsto F_g(p,w)=\exp_{g(p)}w
$$
(outside singular points) where $\exp$ is the exponential map of $\Sf^{2n}$.

As above, any $w\in\Lambda_g$ can be written as
$w=\sum_{j=1}^{n-2}{\mathrm {Re}}(w_jF_j)$
where $w_j=u_j+iv_j$, $1\leq j\leq n-2$.
Then, the submanifold $F_g\colon M^{2n-2}\to\Sf^{2n}$ in \cite{dv2} can be 
parametrized  as 
\bea
F_g(z,w_j)
\!\!\!&=&\!\!\!\frac{1}{\|{\mathrm {Re}}(F_{n+1})\|}
\cos\|\sum_{j=1}^{n-2}(u_j\mathrm{Re}(F_j)
-v_j\mathrm {Im}(F_j))\|\mathrm{Re}(F_{n+1})\\
\!\!\!&&\!\!\!+\,h( \|\sum_{j=1}^{n-2}(u_j\mathrm {Re}(F_j)-v_j\mathrm {Im}(F_j))
\|)\sum_{j=1}^{n-2}(u_j\mathrm {Re}(F_j)-v_j\mathrm {Im}(F_j))
\eea
where $h(x)=\frac{1}{x}\sin x$ if $x\neq 0$ and $h(0)=1$.

\vspace{.5in} {\renewcommand{\baselinestretch}{1}
\hspace*{-20ex}\begin{tabbing} \indent\= IMPA -- Estrada Dona Castorina, 110
\indent\indent\= Univ. of Ioannina -- Math. Dept. \\
\> 22460-320 -- Rio de Janeiro -- Brazil  \>
45110 Ioannina -- Greece \\
\> E-mail: marcos@impa.br \> E-mail: tvlachos@uoi.gr
\end{tabbing}}
\end{document}